\documentclass{amsart}

\usepackage[totalwidth=17cm,totalheight=24cm]{geometry}
\usepackage[T1]{fontenc}
\usepackage[dvips]{graphicx}
\usepackage{epsfig}
\usepackage{amsmath,amssymb,amscd,amsthm}
\usepackage{subfigure}
\usepackage[latin1]{inputenc}
\usepackage{latexsym}
\usepackage{graphics}
\usepackage{colortbl}
\usepackage{color}
\usepackage{ae}

  \def\cprime{$'$}

\begin{document}
\setcounter{secnumdepth}{3}
\setcounter{tocdepth}{2}
\newtheorem{theorem}{Theorem}
\newtheorem{definition}[theorem]{Definition}
\newtheorem{question}[theorem]{Question}
\newtheorem{lemma}[theorem]{Lemma}
\newtheorem{proposition}[theorem]{Proposition}
\newtheorem{corollary}[theorem]{Corollary}

\title{Existence and uniqueness theorem for convex polyhedral metrics on compact surfaces}
\thanks{For the Proceedings of the Conference Efimov-100. The author want to thank the organizers, in particular Sergey Mikhalev
and Idjad Sabitov. The exposition was improved thanks to comments of Alexander A. Borisenko, Ivan Izmestiev and Jean-Marc Schlenker.}
\author{Fran\c{c}ois Fillastre}
\address{Universit\'e de Cergy-Pontoise,  d\'epartement de math\'ematiques, UMR CNRS 8088,
F-95000 Cergy-Pontoise, France}
\email{francois.fillastre@u-cergy.fr}
\date{\today}
\maketitle
\begin{abstract}
 We state that any constant curvature Riemannian metric with conical singularities of constant sign curvature on
 a compact (orientable) surface $S$ can be realized as a convex polyhedron in a 
  (Riemannian or Lorentzian) space-form. Moreover such a polyhedron is unique, up to global isometries, 
among convex polyhedra invariant under isometries acting on a totally umbilical surface.

This general statement falls apart into 10 different cases. 
The cases when $S$ is the sphere are classical.
\end{abstract}

\section{Original statement for the sphere}

Let $P$ be (the boundary of) a convex polytope in the Euclidean space $\mathbb{R}^3$.
The induced metric on $P$ is isometric to a flat metric on the sphere $\mathbb{S}^2$, with conical singularities, 
corresponding to the vertices. As $P$ is convex, the angle around each conical singularity is lower than $2\pi$, and 
then the singular curvature ($2\pi$ minus the cone-angle) is positive. For simplicity we will call such a metric on $\mathbb{S}^2$
a $(0,+)$-metric.
The following famous theorem of A.D.~Alexandrov, proved in the early  1940s, states
the converse property. 
\begin{theorem}[{Alexandrov Theorem, \cite{ale05}}]\label{thm:alex}
 Let $m$ be a  
$(0,+)$-metric on the sphere $\mathbb{S}^2$. Then
\begin{enumerate}
 \item There exists a convex polytope $P$ in $\mathbb{R}^3$ such that $P$ realizes  $m$.
\item $P$ is unique, up to congruences.
\end{enumerate}
\end{theorem}
(``Realizes'' means that the induced metric on $P$ is isometric
to the metric $g$.) 
Of course, the uniqueness of $P$ in the theorem is only among convex polytopes, as it is easy to construct non-convex polytopes
with the same induced metric than a given convex one. The case of degenerated convex polytopes (doubly covered convex polygons)
 is considered in the theorem above.

A convex polytope can be defined in any  simply connected space of constant curvature. Let us
restrict our attention to curvatures $1$ and $-1$.
The original proof of the preceding theorem extends easily to these cases. The induced metric on a convex polytope in  the sphere $\mathbb{S}^3$ (resp. the hyperbolic space 
$\mathbb{H}^3$) is isometric to a $(1,+)$-metric (resp. a $(-1,+)$-metric) on the sphere $\mathbb{S}^2$: 
a spherical (resp. hyperbolic) metric with conical singularities of positive curvature. In the case of spherical metrics 
we also suppose that the number of conical singularities is $>2$, otherwise it is easy to construct examples for which the uniqueness 
part of the statement below is false \cite{ale05}. 
 
\begin{theorem}[{\cite{ale05}}]\label{thm:alex2}
\begin{description} 
 \item[S] Let $m$ be a  
$(1,+)$-metric on the sphere. Then
\begin{enumerate}
 \item There exists a convex polytope $P$ in $\mathbb{S}^3$ such that $P$ realizes  $m$.
\item $P$ is unique, up to congruences.
\end{enumerate}
\item[H]  Let $m$ be a  
$(-1,+)$-metric on the sphere. Then
\begin{enumerate}
 \item There exists a convex polytope $P$ in $\mathbb{H}^3$ such that $P$ realizes  $m$.
\item $P$ is unique, up to congruences.
\end{enumerate}
\end{description}
\end{theorem}

It is easy to visualize convex polytopes in $\mathbb{S}^3$ or $\mathbb{H}^3$.
Actually those spaces can be defined as (pseudo)-spheres in a vector space:
$$\mathbb{S}^3= \{ (x_1,x_2,x_3,x_4)\vert x_1^2+x_2^2+x_3^2+x_4^2=1\},
\mathbb{H}^3= \{ (x_1,x_2,x_3,x_4)\vert x_1^2+x_2^2+x_3^2-x_4^2=-1, x_1>0\}. $$
 A projection from the origin of (a part) of these spaces onto
any affine plane gives an affine representation of these spaces seen as subspaces of the projective space of dimension $3$.
And a convex polytope is a projective object. Of course, the induced metric is not a projective invariant, 
but from this observation some generalizations of the theorems above become natural.

The classical way of proving these statements is to use a \emph{continuity method}, for which existence part follows from the proof
 of uniqueness part using topological arguments \cite{ale05}. The uniqueness part is a 
 refinement  of the famous Cauchy Theorem about rigidity of convex polytopes. The proof of the latter is based on the so-called 
Cauchy Lemma, which is only of topological nature.
The continuity method can be adapted such that an infinitesimal rigidity result suffices, instead of a rigidity result. 
And infinitesimal rigidity is also of projective nature. This has to be understood in the following sense.  
An infinitesimal isometric deformation of a convex polytope $P$ (in any constant curvature pseudo-Riemannian manifold)
can be defined as  the restriction to each faces of $P$ of Killing fields, such that they coincide on edges. 
If there exists a Killing field such that its restriction to $P$ corresponds to the infinitesimal isometric deformation, 
then the deformation is said to be trivial, and if any deformation is trivial, $P$ is infinitesimally rigid.
But if there exists  a map between two pseudo-Riemannian manifolds sending geodesics to geodesics (it is the case in particular
for projective maps), then there exists 
a map sending Killing fields of one manifold to the Killing fields of the other. This result is based on properties
of the connection \cite{kne30}. It was later rediscovered by Volkov \cite{vol74}, and particular cases were known for a long time
(Darboux--Sauer, Pogorelov, $\ldots$).
There is another interpretation of the projective nature of infinitesimal rigidity of convex polytopes, 
based on similar results about frameworks, see \cite{izm09} and references therein.

A proof using the continuity method is not constructive, because it is based on the Domain Invariance Theorem. 
A constructive proof of the existence part of Theorem~\ref{thm:alex} was done in \cite{BI08}. It does not require
a proof of the uniqueness part. It uses a \emph{variational method}, in which the desired polytopes is described as the critical point
of a functional. The definition of this functional uses a polyhedral version of the total scalar curvature. 
Another proof was given in \cite{Vol55}, see the introduction of \cite{BI08}.

\section{Hyperbolic extension}

\begin{figure}[h!]
\begin{center}
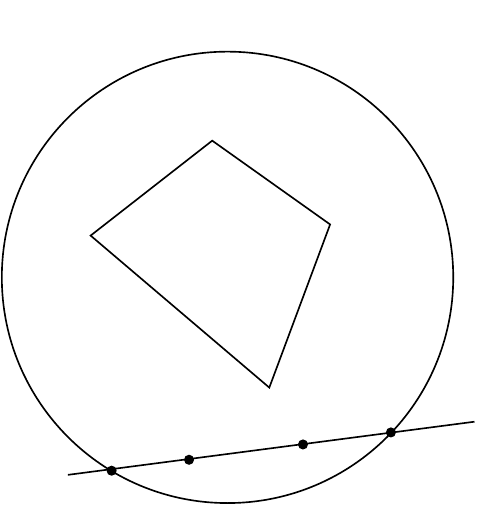 \hspace{2 cm}  
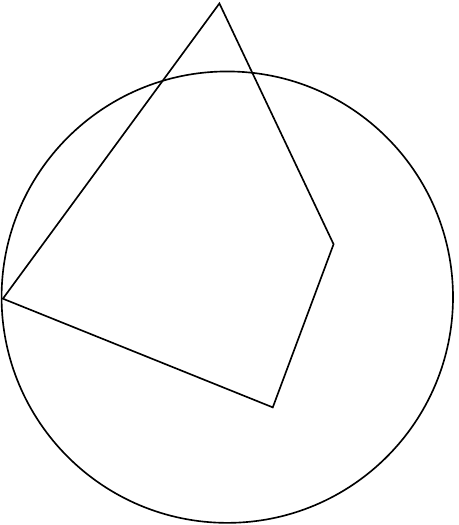 
\end{center}\caption{$d$  is the Hilbert distance, with $[\cdot,\cdot,\cdot,\cdot]$ the cross-ratio.}\label{fig:hyppol}
\end{figure}

The most famous projective model of the hyperbolic space is the Klein model: 
the interior of the unit ball endowed with
the Hilbert distance (Figure~\ref{fig:hyppol}, left). Any 
(projective) convex polytope contained in the interior of the unit ball can be considered as
a hyperbolic polytope. But one can look at the induced hyperbolic metric  on the intersection 
of the unit ball with convex polytopes that are not 
 totally contained in the unit ball (Figure~\ref{fig:hyppol}, right). If each edge of the projective polytope meets the interior
of the unit ball, we call its hyperbolic part a  \emph{generalized hyperbolic polyhedron}.
The induced metric on the neighborhood of a vertex on the unit sphere is isometric to a hyperbolic cusp. 
If the vertex is out of the closed unit ball, the induced metric on the intersection of the faces meeting at this vertex with
the hyperbolic space is isometric to a complete end of infinite area.

Hyperbolic part of Theorem~\ref{thm:alex2} can be extended in this way:

\begin{theorem}\label{ref:sphere hyp gen}
 Let $m$ be a  hyperbolic metric on 
 the (punctured) sphere $\mathbb{S}^2$ with
conical singularities of positive singular curvature, cusps and complete ends of infinite area. Then
\begin{enumerate}
 \item There exists a convex generalized hyperbolic polyhedron $P$ in $\mathbb{H}^3$ such that $P$ realizes  $m$.
\item $P$ is unique, up to congruences.
\end{enumerate}
\end{theorem}

The case with only cusps was proved in \cite{riv94}--- this reference also contains the uniqueness part of
the theorem above.
The proof of the theorem is a part of a more general result contained in \cite{sch98}, see Section~\ref{sect rel}.
A constructive proof for the case with cusps has been announced, see \cite{BPS10}.

\section{De Sitter extension}

\begin{figure}[h!]
\begin{center}
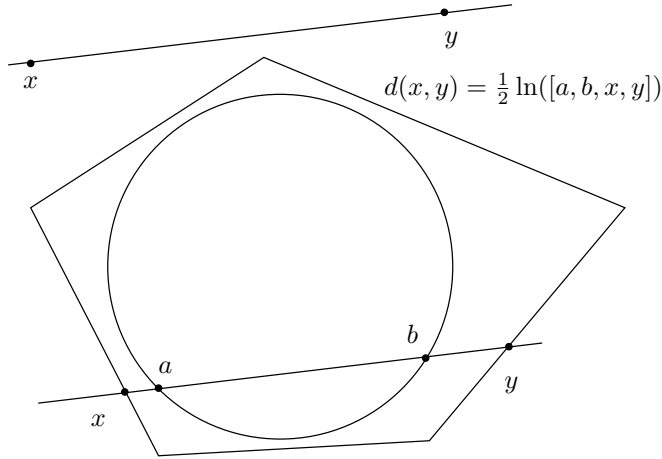 
\caption{When the line don't meet the closed ball, the space is complexified. Then intersection points with the unit sphere exist,
and a imaginary ``distance'' can be defined.}\label{fig de sit}
\end{center}
\end{figure}

The next step is to realize Riemannian metrics in a Lorentzian space.
The outside of the projective model of the hyperbolic space is a projective model of (a part of) de Sitter space, see Figure~\ref{fig de sit}. De Sitter space
is the simply-connected Lorentzian space with constant curvature $1$. It can also be defined as a pseudo-sphere in the Minkowski space,
 endowed with the induced metric:
$$ \mathrm{dS}^3= \{ (x_1,x_2,x_3,x_4)\vert x_1^2+x_2^2+x_3^2-x_4^2=1\}.$$

Consider a convex polytope $P$ in $\mathbb{H}^3$. In the projective Klein model, $P$ has a polar dual
with respect to the unit sphere, denoted by $P^*$. We can see $P^*$ as a convex polyhedral surface in de Sitter space.
The induced metric on $P^*$ is isometric to a spherical (Riemannian) metric on the sphere $\mathbb{S}^2$ with conical
singularities of negative curvature. Moreover the length of closed geodesics is $>2\pi$. We call such a metric a \emph{large}
$(1,-)$-metric on $\mathbb{S}^2$. In particular the number of cone-points is $>2$. Actually the metric induced on the dual $P^*$ of $P$ is isometric to the dual metric, or Gauss image
(the polyhedral analogue of the third fundamental form), of $P$. See \cite{RH93,CD95}.

\begin{theorem}[{\cite{rivinthese,RH93}}]\label{thm:RH}
  Let $m$ be a large $(1,-)$-metric  
 on the sphere. Then
\begin{enumerate}
 \item There exists a convex polyhedral surface $P$ in de Sitter space such that $P$ realizes  $m$.
\item $P$ is unique, up to congruences.
\end{enumerate}
\end{theorem}

\begin{figure}[h!]
\begin{center}
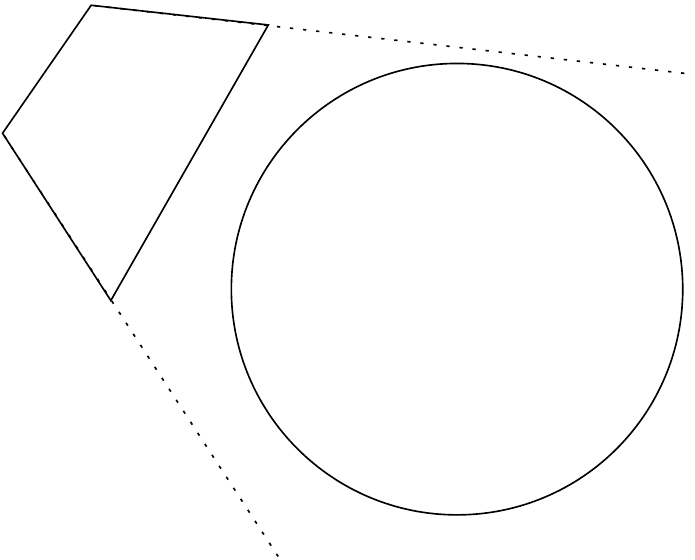 
\caption{A convex polytope in de Sitter space whose induced metric is Riemannian, but which is not the dual of a hyperbolic convex polytope.}
\end{center}\label{fig:desitcpt}
\end{figure}

Actually the uniqueness statement proved in \cite{RH93} is among convex polyhedral surface of de Sitter space which are polar dual 
of hyperbolic convex polytopes. The uniqueness among all convex polyhedral surfaces (homeomorphic to the sphere) of de Sitter space
is a consequence of the results of  \cite{sch01} (in other words, only  dual of convex polytopes of 
$\mathbb{H}^3$ have induced large $(1,-)$-metric). This is not straightforward as there exists convex polytopes in de Sitter space (in the sense that
they bound a compact set of de Sitter space) whose induced metric is Riemannian, see Figure~\ref{fig:desitcpt}. They are characterized by the 
following theorem.

\begin{theorem}[{\cite{sch01}}]
 Let $m$ be a metric on the sphere $\mathbb{S}^2$ such that:
\begin{itemize}
\renewcommand{\labelitemi}{$\bullet$}
 \item $m$ is a spherical metric except at a finite number of singular points $x_1,\dots,x_N$;
\item $\mathbb{S}^2=D_+\cup D_-$, where $D_+$ and $D_-$ are topological discs intersecting on a curve
$C=\partial D_+=\partial D_-$ of length $<2\pi$ with vertices the singular points $x_1,\ldots,x_p$;
\item the singular curvature of $m$ at $x_{p+1},\ldots,x_N$ is negative;
\item $D_+, D_-$ are convex for $m$, and strictly convex at each vertex of $C$;
\item geodesics segments of $D_+,D_-$ have length less than $\pi$. 
\end{itemize}
Then
\begin{enumerate}
 \item There exists a convex polytope $P$ in de Sitter space  such that $P$ realizes  $m$.
\item $P$ is unique, up to congruences.
\end{enumerate}
\end{theorem}

A similar statement holds in the Minkowski space $\mathbb{R}^{2,1}$.

\begin{theorem}[{\cite{IS90,sch01}}]
  Let $m$ be a metric on the sphere $\mathbb{S}^2$ such that:
\begin{itemize}
\renewcommand{\labelitemi}{$\bullet$}
 \item $m$ is a flat metric except at a finite number of singular points $x_1,\dots,x_N$;
\item $\mathbb{S}^2=D_+\cup D_-$, where $D_+$ and $D_-$ are topological discs intersecting on a curve
$C=\partial D_+=\partial D_-$  with vertices the singular points $x_1,\ldots,x_p$;
\item the singular curvature of $m$ at $x_{p+1},\ldots,x_N$ is negative;
\item $D_+, D_-$ are convex for $m$, and strictly convex at each vertex of $C$.
\end{itemize}
Then
\begin{enumerate}
 \item There exists a convex polytope $P$ in $\mathbb{R}^{2,1}$ such that $P$ realizes  $m$.
\item $P$ is unique, up to congruences.
\end{enumerate}
\end{theorem}

As far I know, there does not exist a similar statement for convex polytopes with induced Riemannian metric in the \emph{anti-de Sitter space},
which is a Lorentzian space-form (non-simply connected) of curvature $-1$, which can  be defined as the following pseudo-sphere, endowed 
with the induced metric:
$$ \mathrm{AdS}^3= \{ (x_1,x_2,x_3,x_4)\vert x_1^2+x_2^2-x_3^2-x_4^2=-1\}.$$
Note that the two last theorems are somewhat different of the preceding ones, 
as they don't describe the induced metric on all
the convex polytopes (this one can be non Riemannian for a random convex polytope). Moreover 
the induced metrics have both positive and negative singular curvatures. Anyway
 they are specific to the case of the
sphere.

\section{Higher genus?}

The previous theorems (\ref{thm:alex}, \ref{thm:alex2}, \ref{thm:RH}) are concerned with some 
 $(K,\epsilon)$-metric on the sphere, $K\in\{-1,0,1\}, 
\epsilon\in\{-,+\}$. For which $K$ and $\epsilon$ such metrics exist? The answer is given by Gauss--Bonnet formula
\begin{equation}\label{GB}
 \sum_i k_i=2\pi(2-2g) - KA(S),
\end{equation}
where $g$ is the genus of the compact surface $S$ ($g=0$ in the present case $S$ is the sphere), $k_i$ are the singular curvatures
and $A(S)$ is the area of $S$. Hence there is no other  $(K,\epsilon)$-metric on the sphere than the ones previously considered:
 $(0,+), (1,+), (-1,+), (1,-)$.

So a natural question is: what about realization of $(K,\epsilon)$-metrics on compact surfaces of higher genus?
The aim is to realize such metrics as convex polyhedra. Hence a  $(K,\epsilon)$-metric has to be realized in
$M_K^{\epsilon}$, where
$$M_0^+=\mathbb{R}^3, M_1^+=\mathbb{S}^3, M_{-1}^+=\mathbb{H}^3, M_0^-=\mathbb{R}^{2,1}, M_1^-=\mathrm{dS}^3, M_{-1}^-=\mathrm{AdS}^3.$$
But each of these spaces (or parts of these spaces) has a model in an affine chart of the projective space, and in 
this space convexity and  compactness implies genus $0$.
Then we will realize the universal cover of the metric. We will obtain
convex polyhedral surface equivariant under the action of a representation of the fundamental group of the surface in the isometry group
of (the suitable) $M_K^{\epsilon}$.

We give some examples before definitions and   statements.

\section{Example of flat surfaces of genus $2$}

The most simple example of compact hyperbolic surface is given by the quotient of $\mathbb{H}^2$
by a group of isometries  whose a fundamental domain is a regular octagon (in the Klein projective model).
Such a group is a Fuchsian group, that explains the terminology used below.
This group has a representation $\Gamma$ into the group of the isometries of Minkowski space preserving the
upper-sheet of the hyperboloid. 

Let $p\in\mathbb{H}^2$ seen as a vector of the Minkowski space.
We call $P$  the convex hull \emph{in the ambient Minkowski space} of the orbit 
of $p$ under the action of $\Gamma$ (that gives an infinite number of points on the hyperboloid).

\begin{center}
\begin{figure}
\includegraphics[height=4cm]{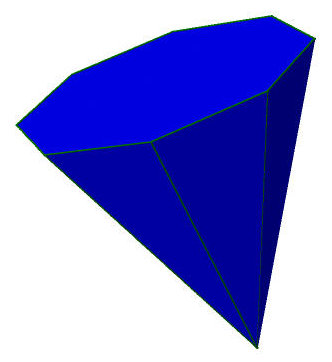} 
 \includegraphics[height=6cm]{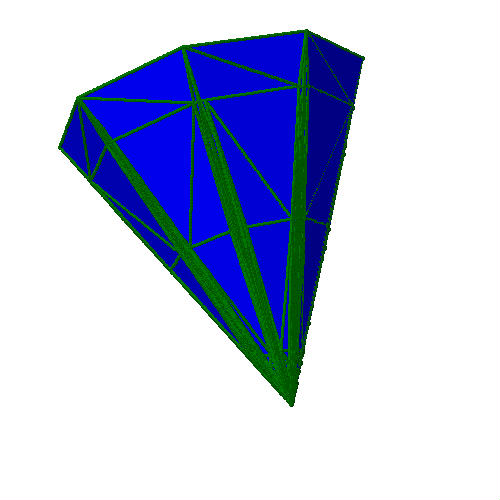}
\includegraphics[height=3.5cm]{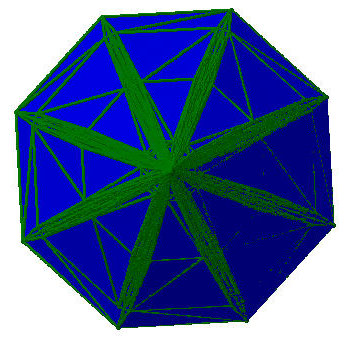}
\includegraphics[height=3.5cm]{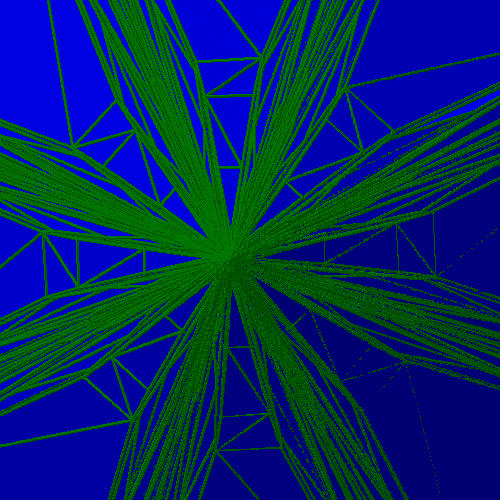}
\includegraphics[height=3.5cm]{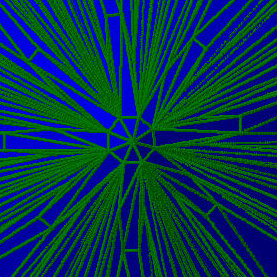}
\caption{The first picture is the convex hull of $p=(0,0,1)$ together with its images for the generators of $\Gamma$
(four elements together with their inverses). The second picture is the convex hull of $p$ and points 
$\gamma_1\circ \gamma_2\circ \gamma_3 (p)$, where $\gamma_i$ is a generator of $\Gamma$ (or the inverse of a generator).
The height (the third coordinate) of $p$ is $1$. The height of $\gamma_1 (p)$ is $\approx 10$. The height of
$\gamma_1 \circ \gamma_2 (p)$ is $\approx 100$. The height  of  $\gamma_1\circ \gamma_2\circ \gamma_3 (p)$
is $\approx 1000$. The remaining pictures are zoomed views of the surface from the bottom.
 For precision reasons the vertices in the images above belongs to the hyperboloid of ray $1000$.
These objects can be manipulate at the author's web page. They were realized with Sage \cite{sage}.}\label{fig:polfuch}
\end{figure}
\end{center}

By construction, $P$ is a convex polyhedral surface, invariant under the action of $\Gamma$, and 
it can be shown that the induced metric is Riemannian. The pair $(P,\Gamma)$ is an example of 
\emph{Fuchsian polyhedron}.
The quotient $P/\Gamma$ is isometric to a compact surface of genus $2$ with a flat metric with one conical singularity 
of curvature $2\pi-8\times \frac{3\pi}{4}=-4\pi<0$ (cone-angle $6\pi$). 
Parts of such a polyhedral surface are drawn on Figure~\ref{fig:polfuch}, with $p=(0,0,1)$.
This is an example of the so-called \emph{convex hull construction}, 
which appeared in the seminal papers \cite{pen87,EP88,NP91}.
The definition of convex Fuchsian polyhedron in the Minkowski space is as follows. The vertices are orbits of points which can live
on different upper-sheets of hyperboloids centered at the origin.
 \begin{definition}
 A \emph{convex Fuchsian polyhedron} in the Minkowski space is a pair $(P,\Gamma)$, where
\begin{itemize}
\renewcommand{\labelitemi}{$\bullet$}
 \item $P$ is a space-like (i.e.~Riemannian) convex polyhedral surface in the Minkowski space $\mathbb{R}^{2,1}$,
\item $\Gamma$ is a group of isometries of $\mathbb{R}^{2,1}$ acting cocompactly on $\mathbb{H}^2\subset \mathbb{R}^{2,1}$
(i.e.~$\mathbb{H}^2/\Gamma$ is a compact hyperbolic surface),
\end{itemize}
 such that $\Gamma$ leaves $P$ (globally) invariant. 
 \end{definition}
And the corresponding existence and uniqueness statement reads as:
\begin{theorem}[{\cite{fil10}}]
 Let $m$ be a  
$(0,-)$-metric on a compact surface $S$ of genus $>1$. Then
\begin{enumerate}
 \item There exists a convex Fuchsian polyhedron $(P,\Gamma)$ in $\mathbb{R}^{2,1}$ such that $P$ realizes  $m$
(i.e.~$P/\Gamma$ endowed with the induced metric is isometric to $(S,m)$).
\item $(P,\Gamma)$ is unique (among convex Fuchsian polyhedra), up to congruences.
\end{enumerate}
\end{theorem}

\section{Example of hyperbolic torus}

The hyperbolic space has a model (the upper-space model) in which
horizontal planes (in the upper part of space) are isometric to the Euclidean plane. They are examples
of horospheres. Let us consider the plane $H$ at height $1$, and its tessellation by 
squares of unit length. The quotient of the Euclidean plane by the group whose fundamental domains are the squares
is a flat torus. This group acting on the Euclidean plane extends to a group $\Gamma$ of  isometries of the hyperbolic space 
(the group is made of parabolic isometries). There is a map from this model of hyperbolic space to
the Klein projective model, which sends $H$ to an ellipsoid of rays $(1/2,1/4)$, passing through the origin and tangent
to the unit sphere. We also denote by $H$ this ellipsoid. Let $p\in H$. The orbit 
$\Gamma p$ is made of an infinite number of points accumulating on the point of tangency between $H$ and the unit sphere.

\begin{figure}
\begin{center}\includegraphics[height=5cm]{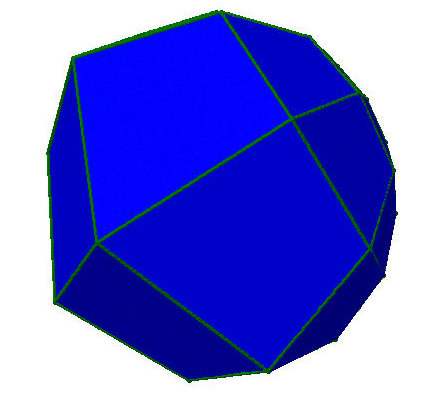}
\includegraphics[height=5cm]{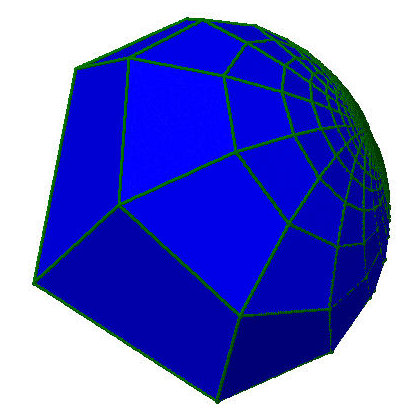}
\includegraphics[height=5cm]{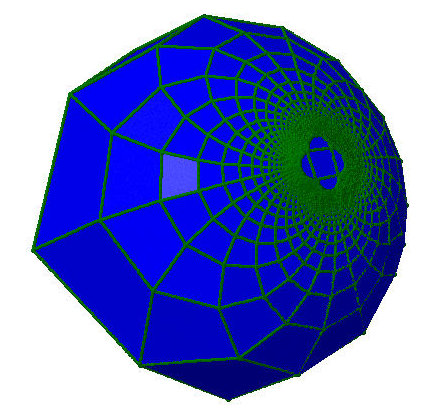}
\includegraphics[height=5cm]{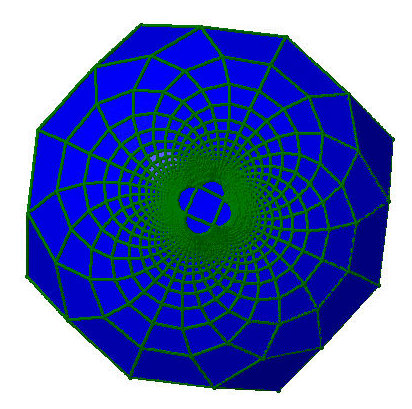}\end{center}
\caption{The parabolic polyhedron is made of an infinite number of vertices accumulating on a same point.
(Only a finite number of vertices are represented on the picture.) 
This object can be manipulate at the author's web page. It was realized with Sage \cite{sage}.}\label{fig:polpar}
\end{figure}

We call $P$  the convex hull in the hyperbolic space (that is the same thing, in the Klein projective model, than
the convex hull in $\mathbb{R}^3$) of $\Gamma p$.
By construction, $P$ is a convex polyhedral surface, invariant under the action of $\Gamma$. All the faces are isometric
hyperbolic squares.
The pair $(P,\Gamma)$ is an 
example of 
\emph{parabolic polyhedron} of the hyperbolic space. 
The quotient $P/\Gamma$ is isometric to a torus with a hyperbolic metric with one conical singularity of
positive curvature (obtained by identifying opposed edges of any face). 
Such a polyhedral surface is drawn on Figure~\ref{fig:polpar}, with $p$ the origin.
Note that we can associate to the same polyhedral surface different parabolic polyhedra, by changing
the group action, such that the different quotients are not isometric. 
The definition of convex parabolic polyhedron in the hyperbolic space is as follows. 
The vertices are orbits of points which can live
on different ellipsoids having the same point of tangency with the unit sphere.
 \begin{definition}
 A \emph{convex parabolic polyhedron} in the 
hyperbolic space is a pair $(P,\Gamma)$, where
\begin{itemize}
\renewcommand{\labelitemi}{$\bullet$}
 \item $P$ is a convex polyhedral surface in $\mathbb{H}^{3}$,
\item $\Gamma$ is a group of isometries of $\mathbb{H}^{3}$ acting cocompactly on a horosphere,
\end{itemize}
 such that $\Gamma$ leaves $P$ (globally) invariant. 
 \end{definition}
And the corresponding existence and uniqueness statement reads as:
\begin{theorem}[{\cite{FI09}}]\label{thm:par}
 Let $m$ be a  
$(-1,+)$-metric on the torus $T$. Then
\begin{enumerate}
 \item There exists a convex parabolic polyhedron $(P,\Gamma)$ in $\mathbb{H}^{3}$ such that $P$ realizes  $m$
(i.e.~$P/\Gamma$ endowed with the induced metric is isometric to $(T,m)$).
\item $(P,\Gamma)$ is unique (among convex parabolic polyhedra), up to congruences.
\end{enumerate}
\end{theorem}

The uniqueness have to be among convex parabolic polyhedra, as there is another way to realize 
$(-1,+)$-metrics on the torus. Let us call \emph{banana} a surface of $\mathbb{H}^3$ at constant distance
from a geodesic. The banana endowed with the induced metric is isometric 
to a complete infinite Euclidean cylinder.  It is then easy to construct 
a polyhedral surface $P$ and a group of isometries $\Gamma$ acting cocompactly on the banana, such that $P/\Gamma$ is
isometric to a $(-1,+)$-metric on the torus. If we call such a pair a \emph{convex polyhedral banana},
the following question holds
\begin{question}\label{ques bana}
 Let $m$ be a  
$(-1,+)$-metric on the torus $T$. Does it exist a (unique) convex polyhedral banana realizing $(T,m)$?
\end{question}
It seems that there is no other way to obtain $(-1,+)$-metrics on the torus, as
horospheres and bananas are the only complete isometric immersions of the Euclidean plane
into $\mathbb{H}^3$ \cite{VV71,sas73}. The same remarks can be done  in the case of de Sitter space \cite{DN81}, in which horospheres and bananas 
can be defined similarly.

\section{How many cases?}

What are the possible $(K,\epsilon)$-metric 
on a surface of genus $g$? A priori there are $18$ cases:
$3$ choices for $K$, $2$ choices for $\epsilon$, and $3$ choices for
$g$ (depending if the universal cover is a sphere, a disc or a plane).
But by Gauss--Bonnet Formula (\ref{GB}), there are only $10$ possible cases, which are listed in Table~\ref{tableau}.

\renewcommand{\arraystretch}{1.3}
\begin{center}\begin{table}[h!]
\begin{tabular}{| c | c | c |  c | c |}
   \hline  $g$ & $K$ & $\epsilon$ &  DM & VM \\
 \hline    $0$ & $-1$ & $+$  & \cite{alex42} & \cite{BI08} \\
   $0$ & $0$ & $+$  &\cite{ale05}  & \\
     $0$ & $1$ & $+$  &\cite{ale05}  & \\
    $0$ & $1$ & $-$   &\cite{RH93} &\\
   \hline    $1$ & $-1$ & $+$ &  &\cite{FI09}\\
    $1$ & $1$ & $-$ & &  \cite{FI10} \\
   \hline   $\geq 2$& $-1$ & $+$  &\cite{artrealisationhyp} &\\
    $\geq 2$ & $-1$& $-$ & \cite{fil10}& \\
     $\geq 2$ & $0$& $-$& \cite{fil10} &\\
      $\geq 2$  & $1$ & $-$   & \cite{sch04,fil10} &\\ \hline 
      \end{tabular}
 \caption{\label{tableau} The ten cases of realization of metrics of curvature $K$ with cone singularities of curvature of 
sign $\epsilon$ on a compact surface of genus $g$. The column named by DM lists the articles using a deformation method and the column 
named by VM lists the articles using a variational method. The $(1,-)$-metrics are supposed to be large.}
\end{table}\end{center}

\section{The general statement}

The ten statements listed in Table~\ref{tableau} can be resumed in a unique statement.
Let $S$ be a compact surface of genus $g$,  $K\in\{-1,0,1\}$ and $\epsilon\in\{-,+\}$.
Sometimes $(0,+)$-metrics on the sphere were called \emph{convex polyhedral metrics} on the sphere, as
by Alexandrov Theorem they correspond to induced metrics on convex polytopes in Euclidean space. By analogy we call
any $(K,\epsilon)$-metric a convex polyhedral metric, as the following statement says that they correspond 
to induced metrics on convex polyhedral surfaces.

\begin{theorem}[Existence and uniqueness theorem for convex polyhedral metrics on compact surfaces]
 Let $m$  be a $(K,\epsilon)$-metric on $S$ (large for $(K,\epsilon=(1,-)$).
 Then
\begin{enumerate}
 \item There exists a convex polyhedral surface $P\subset M_K^{\epsilon}$
and a group $\Gamma$ of isometries $\mathbb{R}^{2,1}$, acting cocompactly on a totally umbilic surface, such that $P/\Gamma$ is isometric to 
$(S,m)$.
\item $(P,\Gamma)$ is the unique pair of this kind, up to congruences.
\end{enumerate}
\end{theorem}

If, for a pair $(P,\Gamma)$ as in the statement above, the umbilic surface is isometric to the Euclidean plane,
$(P,\Gamma)$ is called a convex parabolic polyhedron, and if the umbilic surface is isometric to the hyperbolic plane (up to a homothety), 
$(P,\Gamma)$ is called a convex Fuchsian polyhedron. If $S$ is the sphere, the action of $\Gamma$ is trivial.

In the $(-1,+)$-cases, the statement above can be extended in the same way than Theorem~\ref{ref:sphere hyp gen}.
\begin{theorem}[{\cite{fil08}}]\label{thm:gen hyp}
 Let $m$ be a  hyperbolic metric on 
 a surface $S$ of finite topological type with
conical singularities of positive singular curvature, cusps and complete ends of infinite area. Then
\begin{enumerate}
 \item There exists a convex generalized hyperbolic polyhedral surface $P$ in $\mathbb{H}^3$ 
and a group $\Gamma$ of isometries of $\mathbb{H}^3$, acting cocompactly
on a totally umbilic surface, such that $P/\Gamma$ is isometric to $(S,m)$.
\item $(P,\Gamma)$ is the unique couple of this kind, up to congruences.
\end{enumerate}
\end{theorem}

The uniqueness in the statement above is only among polyhedral surfaces invariant under the 
action of these groups acting on umbilic surfaces. In a certain sense they are the ``most simple''. For example
it is possible to construct convex polyhedral surfaces in $\mathbb{H}^3$ invariant under groups more complicated than 
Fuchsian groups (as for example quasi-Fuchsian groups), whose quotient provides $(-1,+)$-metrics on compact surfaces of genus $>1$.

\section{Related works}\label{sect rel}

\paragraph{\textbf{Higher dimensions}} The results presented here are specific to surfaces. A result such as 
Theorem~\ref{thm:alex} is possible because on one hand the space of convex polytopes with 
$n$ vertices (up to congruences) has dimension $3n-6$ (the polytope is determined by the positions of the vertices in space,
and isometry group has dimension $6$), and on the other hand the space of 
Euclidean metrics on the sphere with $n$ conical singularities of positive curvature
up to isometries
has dimension $3n-6$ (roughly speaking, the metric is determined by the conformal structure of the punctured sphere
and the values of the singular curvatures \cite{tro86,tro91}). 
Let us discuss for example the Euclidean case of dimension $4$. Convex polytope are still determined by their vertices, but the induced metric
is singular along edges, so the space of polytopes and the space of metrics don't have the same dimension.
Moreover the uniqueness is obvious, at least locally (the link of a vertex is a convex  polytope in $\mathbb{S}^3$, which is rigid
by Theorem~\ref{thm:alex2}, that implies that the Euclidean convex polytope is rigid).

\paragraph{\textbf{Pseudo-Riemannian metrics on surfaces}}
In Theorem~\ref{thm:gen hyp} we considered ``generalized hyperbolic polyhedron'', which are the intersection
of convex polyhedral surface of projective space with the unit ball, endowed with the Hilbert metric. 
But we could look at the induced metric on convex polyhedral surface in the ``hyperbolic-de Sitter'' space
(an extension of $\mathbb{R}^3$ endowed with the Hilbert metric).  
Such induced metric can be
Riemannian, Lorentzian or degenerated on different faces and edges. For closed
polyhedral surfaces with a finite number of vertices, this is done in \cite{sch98,sch01}.
The case of convex polytopes in Minkowski space is done in \cite{sch01}. 

\paragraph{\textbf{Circle patterns}}
The problems of circle patterns on surfaces have deep connections with the geometry of hyperbolic polyhedra. Let us consider the example of the flat torus.
Roughly speaking, if a triangulation $T$ of the torus together with numbers associated to the edges are given, satisfying some 
conditions, then there is a flat metric on the torus (unique up to scalar multiple), a geodesic triangulation isotopic to $T$
and a (unique) family of circles, centered at the vertices of $T$, and intersecting at an angle given by the given numbers 
(this is a \emph{circle pattern} on the torus).
The relation between circle patterns and hyperbolic polyhedra is as follow. Consider a circle pattern of the torus is given, 
and look at the universal cover. 
Let us consider him
as the horizontal plane at height zero in three-space, and draw above each circle a half-sphere. 
In the upper-space model of the hyperbolic space, those half-circles bound
a convex parabolic polyhedron, and the angles between the circles correspond to dihedral angles of the polyhedron.
The triangulation of the circle pattern is the dual of the cellulation of the polyhedron, hence it is simple, and then dihedral angles determine the dual metric:
circle pattern problem is reduced to a realization problem. For example Theorem~\ref{thm:RH} can be considered as a generalization of 
the Koebe--Andreev--Thurston theorem for circle pattern on the sphere (recall that realization of hyperbolic polyhedra with respect to the dual metric
is the same as realization of de Sitter polyhedra with respect to the induced metric).
There exists a variational proof for existence of circle pattern on compact surface \cite{Thurnotes,CL03} but the  functional 
 was not explicitly known.
Actually (in the torus case)
 this functional is the one appearing in the variational proof for realization theorem with respect to
the dual metric \cite{FI10}. It is related to the discrete Hilbert-Einstein functional and uses the volume of hyperbolic polyhedra.

The vertices of the hyperbolic polyhedra constructed from a circle pattern on a compact surface can be inside, on the boundary at infinity, or
outside the hyperbolic space, depending on how the circles meet. Close to the scope of this text are the results of
\cite{BB02,rou04},  about circle patterns corresponding to hyperbolic polyhedra with vertices outside of  hyperbolic space, but with edges meeting
hyperbolic space
(hyperideal polyhedra).

Circle patterns can be generalized to \emph{inversive distance circle packings}, where the inversive distance is a quantity
prescribed between two circles with maybe empty intersection. It the intersection is non-empty, the inversive distance 
is the intersection angle between the circles, see \cite{luo10} and references therein. If the circles are seen on the boundary 
at infinity of the hyperbolic space, the inversive distance between two non-intersecting circles is nothing but the hyperbolic distance 
between the two hyperbolic planes defined by the circles. It would be interesting to know if the results of \cite{luo10}
could be expressed in terms of hyperbolic polyhedra (which are now allowed to have edges outside of hyperbolic space),
 and if the functionals used in this paper could be related to hyperbolic volume.

Connections between realization of ideal hyperbolic polyhedra and ``discrete conformal maps'' are done in \cite{BPS10}.

\paragraph{\textbf{Convex metrics on surfaces}}

Let $m$ be the $(0,+)$-metric induced on a convex polytope of the Euclidean space. 
By the uniqueness part of Theorem~\ref{thm:alex}, there does not exist other convex polytope (up to congruence)
with induced metric $m$. But we can ask if there exists a convex surface with induced metric $m$ 
(a convex surface in $\mathbb{R}^3$ is the boundary of a convex body). The answer is negative, as it was proved by Olovyanishnikov
few times after Alexandrov proved Theorem~\ref{thm:alex}.
This result was extended by Pogorelov: a convex surface is uniquely determined by its induced metric \cite{pog73,bus58}.
 
The induced metric on a convex surface  satisfies
some convexity conditions. A.D.~Alexandrov showed that any metric on the sphere satisfying such conditions
is realized by a convex surface. The proof is obtained from Theorem~\ref{thm:alex}.
by approximation by polytopes. See \cite{ZA67,ale06} for more details. There does not exist such general results in other spaces
as far I know.

Convex polyhedra can be seen as a particular case of convex surfaces. Another particular case are convex smooth surfaces, for which many
results exist. The existence part of Theorem~\ref{thm:alex} is known as Weyl Problem and was solved by many authors, until papers of 
Nirenberg and Pogorelov in the 50's. A variational
proof was proposed in \cite{BH37}, and Ivan Izmestiev is currently working on it. The smooth analogue of
Theorem~\ref{thm:RH} was done in \cite{sch94,sch96}. \cite{sch98b} contains the smooth analogue of Theorem~\ref{ref:sphere hyp gen}.
The smooth analogue of Fuchsian results in the hyperbolic space is a consequence of the more general result of \cite{sch06} (see below).
Fuchsian results in Lorentzian spaces were done in \cite{LS00}. The idea of such results seems to go back to 
\cite[3.2.4,B',B'']{gro86}. Smooth parabolic results are not known for the moment, but a positive answer to the smooth version of
Question~\ref{ques bana} is done in \cite{sch06}. Note that here ``smooth'' means $C^{\infty}$. For a review of results with assumptions
about the regularity of the metric or of the immersion, see \cite{KS91,KS95}.

\paragraph{\textbf{Hyperbolic manifolds with convex boundary}}

The hyperbolic part of Theorem~\ref{thm:alex2}
is equivalent to:
\begin{theorem}
 Let $m$ be a  
$(-1,+)$-metric on the sphere. Then
there exists a unique 
hyperbolic metric on the ball $B$ such that the induced metric
on $\partial B$ is isometric to $m$.
\end{theorem}

Theorem~\ref{thm:RH} provides an analogous statement, with the induced Gauss metric
instead of the induced metric on the boundary.  Theorem~\ref{thm:par}
is equivalent to:
\begin{theorem}
Let $m$ be a  
$(-1,+)$-metric on the torus. Then
there exists a unique 
hyperbolic metric on a cusp $C$ with convex polyhedral boundary such that the induced metric
on $\partial C$ is isometric to $m$.
\end{theorem}
The de Sitter analogue of Theorem~\ref{thm:par} provides a statement
analogous to
 the one above \cite{FI10}.   A natural extension of such statements is to ask if, 
for a hyperbolic 3-manifold with convex boundary, the induced metric on the boundary uniquely
determines the hyperbolic metric. Here there is no other assumption on the boundary than convexity.
However a reasonable assumption on $M$ should be that it admits a complete hyperbolic convex cocompact metric, i.e.~$M$ contains 
a non-empty, compact, geodesically convex subset. This can be though as an extension of Mostow Rigidity Theorem
(the topology and the induced metric on the boundary uniquely determines
the hyperbolic metric on the manifold).  A celebrated particular case is the one of induced metric on the boundary of
the convex core of hyperbolic manifolds. In this case the metric on the boundary is hyperbolic. See the introduction of \cite{sch06} for references.
Besides if the boundary is asked to be smooth and strictly convex, such a statement  was done in 
 \cite{sch06}, answering a question of Thurston. The existence part was proved before in \cite{lab92}. 
Partial results for higher dimension manifolds with smooth boundary are in \cite{sch01b}.

By contrast to the statements of the two theorems  above, which concern topologically simple manifolds, it seems
to be not possible to make a straightforward polyhedral analogue of the smooth statement of \cite{sch06}, because some convex polyhedral metrics  
could be realized by boundaries of  3-manifolds containing some ``pleated'' faces (instead of totally geodesic ones), see \cite{sch03,sch04}.

For each of the problems mentioned above, a ``dual'' problem can be considered as well. For example in the case of the boundary of the convex core, 
the dual statement implies bending laminations on surfaces.

Such considerations can be translated in terms of  ``Globally
hyperbolic AdS manifolds'' (a king of analogue of quasi-Fuchsian manifolds modeled on anti-de Sitter geometry).
In this setting, characterization of the dual metric on the boundary of the convex core is done in \cite{BS10}, and characterization
of the induced metric is the subject of \cite{dial}. 
For more details about this subject, we send the reader to
 the review paper \cite{sch03} and to \cite{mes07,mess-notes}.

Recently,  the kind of statements mentioned here  
are extended from the setting of manifolds to the one of ``manifolds with particles''
(cone-manifolds with singularities along certain geodesics, in particular time-like in the Lorentzian setting), see 
 \cite{BS09,LS09}. Such objects appear naturally in the variational proofs for realization theorems, and 
they are also related to circle patterns problems on surfaces with conical singularities.

\bibliographystyle{alpha}

\begin{thebibliography}{Sch01b}

\bibitem[ABB{\etalchar{+}}07]{mess-notes}
L.~Andersson, T.~Barbot, R.~Benedetti, F.~Bonsante,
  W.~M.~Goldman, F.~Labourie, K.~Scannell, and J.-M.~Schlenker.
\newblock Notes on: ``{L}orentz spacetimes of constant curvature'' [{G}eom.
  {D}edicata {\bf 126} (2007), 3--45; mr2328921] by {G}. {M}ess.
\newblock {\em Geom. Dedicata}, 126:47--70, 2007.

\bibitem[Ale42]{alex42}
A.~D. Alexandroff.
\newblock Existence of a convex polyhedron and of a convex surface with a given
  metric.
\newblock {\em Rec. Math. [Mat. Sbornik] N.S.}, 11(53):15--65, 1942.

\bibitem[Ale05]{ale05}
A.~D. Alexandrov.
\newblock {\em Convex polyhedra}.
\newblock Springer Monographs in Mathematics. Springer-Verlag, Berlin, 2005.
\newblock Translated from the 1950 Russian edition by N. S. Dairbekov, S. S.
  Kutateladze and A. B. Sossinsky, With comments and bibliography by V. A.
  Zalgaller and appendices by L. A. Shor and Yu. A. Volkov.

\bibitem[Ale06]{ale06}
A.~D. Alexandrov.
\newblock {\em A. {D}. {A}lexandrov selected works. {P}art {II}}.
\newblock Chapman \& Hall/CRC, Boca Raton, FL, 2006.
\newblock Intrinsic geometry of convex surfaces, Edited by S. S. Kutateladze,
  Translated from the Russian by S. Vakhrameyev.

\bibitem[AZ67]{ZA67}
A.~D. Aleksandrov and V.~A. Zalgaller.
\newblock {\em Intrinsic geometry of surfaces}.
\newblock Translated from the Russian by J. M. Danskin. Translations of
  Mathematical Monographs, Vol. 15. American Mathematical Society, Providence,
  R.I., 1967.

\bibitem[BB02]{BB02}
X.~Bao and F.~Bonahon.
\newblock Hyperideal polyhedra in hyperbolic 3-space.
\newblock {\em Bull. Soc. Math. France}, 130(3):457--491, 2002.

\bibitem[BH37]{BH37}
W.~Blaschke and G.~Herglotz.
\newblock \"uber die verwirklichung einer geschlossenen fl\"ache mit
  vorgeschriebenem bogenelement im euklidischen raum.
\newblock {\em Sitzungsber. Bayer. Akad. Wiss., Math.-Naturwiss. Abt.},
  2:229--230, 1937.

\bibitem[BI08]{BI08}
A.~I. Bobenko and I.~Izmestiev.
\newblock Alexandrov's theorem, weighted {D}elaunay triangulations, and mixed
  volumes.
\newblock {\em Ann. Inst. Fourier (Grenoble)}, 58(2):447--505, 2008.

\bibitem[BPS10]{BPS10}
A.~Bobenko, U.~Pinkall, and B.~Springborn.
\newblock Discrete conformal maps and ideal hyperbolic polyhedra.
\newblock arXiv:1005.2698, 2010.

\bibitem[BS09]{BS09}
F.~Bonsante and J.-M. Schlenker.
\newblock Ad{S} manifolds with particles and earthquakes on singular surfaces.
\newblock {\em Geom. Funct. Anal.}, 19(1):41--82, 2009.

\bibitem[BS10]{BS10}
F.~Bonsante and J.-M. Schlenker.
\newblock Fixed points of compositions of earthquakes.
\newblock arXiv:0812.3471, 2010.

\bibitem[Bus08]{bus58}
H.~Busemann.
\newblock {\em Convex surfaces}.
\newblock Dover Publications Inc., Mineola, NY, 2008.
\newblock Reprint of the 1958 original.

\bibitem[CD95]{CD95}
R.~Charney and M.~Davis.
\newblock The polar dual of a convex polyhedral set in hyperbolic space.
\newblock {\em Michigan Math. J.}, 42(3):479--510, 1995.
\newblock Correction {\it Michigan Math. J.}, 43(3):619, 1996.

\bibitem[CL03]{CL03}
B.~Chow and F.~Luo.
\newblock Combinatorial {R}icci flows on surfaces.
\newblock {\em J. Differential Geom.}, 63(1):97--129, 2003.

\bibitem[Diaon]{dial}
B.~Diallo.
\newblock {\em Prescription de m\'etriques sur le bords des coeurs convexes de
  vari\'et\'es anti-de {S}itter globalement hyperboliques maximales {C}auchy
  compactes}.
\newblock PhD thesis, Universit\'e de Toulouse, In preparation.

\bibitem[DN81]{DN81}
M.~Dajczer and K.~Nomizu.
\newblock On flat surfaces in {$S^{3}_{1}$} and {$H^{3}_{1}$}.
\newblock In {\em Manifolds and {L}ie groups ({N}otre {D}ame, {I}nd., 1980)},
  volume~14 of {\em Progr. Math.}, pages 71--108. Birkh\"auser Boston, Mass.,
  1981.

\bibitem[EP88]{EP88}
D.~B.~A. Epstein and R.~C. Penner.
\newblock Euclidean decompositions of noncompact hyperbolic manifolds.
\newblock {\em J. Differential Geom.}, 27(1):67--80, 1988.

\bibitem[FI09]{FI09}
F.~Fillastre and I.~Izmestiev.
\newblock Hyperbolic cusps with convex polyhedral boundary.
\newblock {\em Geom. Topol.}, 13(1):457--492, 2009.

\bibitem[FI10]{FI10}
F.~Fillastre and I.~Izmestiev.
\newblock Gauss images of convex polyhedral cusps.
\newblock To appear Transactions of the AMS, 2010.

\bibitem[Fil07]{artrealisationhyp}
F.~Fillastre.
\newblock Polyhedral realisation of hyperbolic metrics with conical
  singularities on compact surfaces.
\newblock {\em Ann. Inst. Fourier (Grenoble)}, 57(1):163--195, 2007.

\bibitem[Fil08]{fil08}
F.~Fillastre.
\newblock Polyhedral hyperbolic metrics on surfaces.
\newblock {\em Geom. Dedicata}, 134:177--196, 2008.
\newblock Erratum {\it Geom. Dedicata}, 138(1):193--194, 2009.

\bibitem[Fil10]{fil10}
F.~Fillastre.
\newblock Fuchsian polyhedra in {L}orentzian space-forms.
\newblock math.DG/0702532. To appear in {\it Mathematische Annalen}, 2010.

\bibitem[Gro86]{gro86}
M.~Gromov.
\newblock {\em Partial differential relations}, volume~9 of {\em Ergebnisse der
  Mathematik und ihrer Grenzgebiete (3) [Results in Mathematics and Related
  Areas (3)]}.
\newblock Springer-Verlag, Berlin, 1986.

\bibitem[IKS91]{KS91}
I.~Ivanova-Karatopraklieva and I.~Kh. Sabitov.
\newblock Deformation of surfaces. {I}.
\newblock In {\em Problems in geometry, {V}ol.\ 23 ({R}ussian)}, Itogi Nauki i
  Tekhniki, pages 131--184, 187. Akad. Nauk SSSR Vsesoyuz. Inst. Nauchn. i
  Tekhn. Inform., Moscow, 1991.
\newblock Translated in J. Math. Sci. {{\bf{7}}0} (1994), no. 2, 1685--1716.

\bibitem[IKS95]{KS95}
I.~Ivanova-Karatopraklieva and I.~Kh. Sabitov.
\newblock Bending of surfaces. {II}.
\newblock {\em J. Math. Sci.}, 74(3):997--1043, 1995.
\newblock Geometry, 1.

\bibitem[IS90]{IS90}
U.~Il{\cprime}khamov and D.~D. Sokolov.
\newblock Realizability of lentil-shaped polyhedra in a pseudo-{E}uclidean
  space.
\newblock {\em Vestnik Moskov. Univ. Ser. I Mat. Mekh.}, (2):3--6, 104, 1990.

\bibitem[Izm09]{izm09}
I.~Izmestiev.
\newblock Projective background of the infinitesimal rigidity of frameworks.
\newblock {\em Geom. Dedicata}, 140:183--203, 2009.

\bibitem[Kne30]{kne30}
M.~S. Knebelman.
\newblock On {G}roups of {M}otion in {R}elated {S}paces.
\newblock {\em Amer. J. Math.}, 52(2):280--282, 1930.

\bibitem[Lab92]{lab92}
F.~Labourie.
\newblock M\'etriques prescrites sur le bord des vari\'et\'es hyperboliques de
  dimension {$3$}.
\newblock {\em J. Differential Geom.}, 35(3):609--626, 1992.

\bibitem[LS00]{LS00}
F.~Labourie and J.-M. Schlenker.
\newblock Surfaces convexes fuchsiennes dans les espaces lorentziens \`a
  courbure constante.
\newblock {\em Math. Ann.}, 316(3):465--483, 2000.

\bibitem[LS09]{LS09}
C.~Lecuire and J.-M. Schlenker.
\newblock The convex core of quasifuchsian manifolds with particles, 2009.
\newblock arXiv:0909.4182.

\bibitem[Luo10]{luo10}
F.~Luo.
\newblock Rigidity of polyhedral surfaces, {III}.
\newblock arXiv:1010.3284, 2010.

\bibitem[Mes07]{mes07}
G.~Mess.
\newblock Lorentz spacetimes of constant curvature.
\newblock {\em Geom. Dedicata}, 126:3--45, 2007.



\bibitem[NP91]{NP91}
M.~N{\"a}{\"a}t{\"a}nen and R.~C. Penner.
\newblock The convex hull construction for compact surfaces and the {D}irichlet
  polygon.
\newblock {\em Bull. London Math. Soc.}, 23(6):568--574, 1991.

\bibitem[Pen87]{pen87}
R.~C. Penner.
\newblock The decorated {T}eichm\"uller space of punctured surfaces.
\newblock {\em Comm. Math. Phys.}, 113(2):299--339, 1987.

\bibitem[Pog73]{pog73}
A.~V. Pogorelov.
\newblock {\em Extrinsic geometry of convex surfaces}.
\newblock American Mathematical Society, Providence, R.I., 1973.
\newblock Translated from the Russian by Israel Program for Scientific
  Translations, Translations of Mathematical Monographs, Vol. 35.

\bibitem[RH93]{RH93}
I.~Rivin and C.~D. Hodgson.
\newblock A characterization of compact convex polyhedra in hyperbolic
  {$3$}-space.
\newblock {\em Invent. Math.}, 111(1):77--111, 1993.

\bibitem[Riv86]{rivinthese}
I.~Rivin.
\newblock {\em On geometry of convex polyhedra in hyperbolic 3-space}.
\newblock PhD thesis, Princeton University, June 1986.

\bibitem[Riv94]{riv94}
I.~Rivin.
\newblock Intrinsic geometry of convex ideal polyhedra in hyperbolic
  {$3$}-space.
\newblock In {\em Analysis, algebra, and computers in mathematical research
  (Lule\aa, 1992)}, volume 156 of {\em Lecture Notes in Pure and Appl. Math.},
  pages 275--291. Dekker, New York, 1994.

\bibitem[Rou04]{rou04}
M.~Rousset.
\newblock Sur la rigidit\'e de poly\`edres hyperboliques en dimension 3: cas de
  volume fini, cas hyperid\'eal, cas fuchsien.
\newblock {\em Bull. Soc. Math. France}, 132(2):233--261, 2004.

\bibitem[S{\etalchar{+}}10]{sage}
W.\thinspace{}A. Stein et~al.
\newblock {\em {S}age {M}athematics {S}oftware ({V}ersion 4.5.3)}.
\newblock The Sage Development Team, 2010.
\newblock {\tt http://www.sagemath.org}.

\bibitem[Sas73]{sas73}
S.~Sasaki.
\newblock On complete flat surfaces in hyperbolic {$3$}-space.
\newblock {\em K\=odai Math. Sem. Rep.}, 25:449--457, 1973.

\bibitem[Sch94]{sch94}
J.-M. Schlenker.
\newblock Surfaces elliptiques dans des espaces lorentziens \`a courbure
  constante.
\newblock {\em C. R. Acad. Sci. Paris S\'er. I Math.}, 319(6):609--614, 1994.

\bibitem[Sch96]{sch96}
J.-M. Schlenker.
\newblock Surfaces convexes dans des espaces lorentziens \`a courbure
  constante.
\newblock {\em Comm. Anal. Geom.}, 4(1-2):285--331, 1996.

\bibitem[Sch98a]{sch98}
J.-M. Schlenker.
\newblock M\'etriques sur les poly\`edres hyperboliques convexes.
\newblock {\em J. Differential Geom.}, 48(2):323--405, 1998.

\bibitem[Sch98b]{sch98b}
J.-M. Schlenker.
\newblock R\'ealisations de surfaces hyperboliques compl\`etes dans {$H^3$}.
\newblock {\em Ann. Inst. Fourier (Grenoble)}, 48(3):837--860, 1998.

\bibitem[Sch01a]{sch01}
J.-M. Schlenker.
\newblock Convex polyhedra in {L}orentzian space-forms.
\newblock {\em Asian J. Math.}, 5(2):327--363, 2001.

\bibitem[Sch01b]{sch01b}
J.-M. Schlenker.
\newblock Einstein manifolds with convex boundaries.
\newblock {\em Comment. Math. Helv.}, 76(1):1--28, 2001.

\bibitem[Sch03]{sch03}
J.-M. Schlenker.
\newblock Hyperideal polyhedra in hyperbolic manifolds.
\newblock Preprint, 2003.

\bibitem[Sch04]{sch04}
J.-M. Schlenker.
\newblock Hyperbolic manifolds with polyhedral boundary.
\newblock Preprint, 2004.

\bibitem[Sch06]{sch06}
J.-M. Schlenker.
\newblock Hyperbolic manifolds with convex boundary.
\newblock {\em Invent. Math.}, 163(1):109--169, 2006.

\bibitem[Thu97]{Thurnotes}
W.~P. Thurston.
\newblock {\em The geometry and topology of three-manifolds. Recent version of
  the 1980 notes.}
\newblock http://www.msri.org/gt3m, 1997.

\bibitem[Tro86]{tro86}
M.~Troyanov.
\newblock Les surfaces euclidiennes \`a singularit\'es coniques.
\newblock {\em Enseign. Math. (2)}, 32(1-2):79--94, 1986.

\bibitem[Tro91]{tro91}
M.~Troyanov.
\newblock Prescribing curvature on compact surfaces with conical singularities.
\newblock {\em Trans. Amer. Math. Soc.}, 324(2):793--821, 1991.

\bibitem[Vol55]{Vol55}
{Yu.}~A. Volkov.
\newblock {\em Existence of a polyhedron with a given development}.
\newblock PhD thesis, Leningrad State University, 1955.
\newblock Russian.

\bibitem[Vol74]{vol74}
Yu.~A. Volkov.
\newblock A generalization of the {D}arboux-{S}auer and {P}ogorelov theorems.
\newblock {\em Zap. Nau\v cn. Sem. Leningrad. Otdel. Mat. Inst. Steklov.
  (LOMI)}, 45:63--67, 118, 1974.
\newblock Problems in global geometry.

\bibitem[VV71]{VV71}
Yu.~A. Volkov and S.~M. Vladimirova.
\newblock Isometric immersions of a euclidean plane in lobachevskii space.
\newblock {\em Mathematical Notes}, 10(3):619--622, 1971.

\end{thebibliography}
\newcommand{\etalchar}[1]{$^{#1}$}

\end{document}